\documentclass[smallextended]{svjour3}       
\smartqed  
\usepackage{graphicx}
\usepackage{mathrsfs} 
\usepackage{amsmath} 
\usepackage{csquotes} 
\usepackage{mathtools} 
\usepackage{amssymb}
\usepackage{tikz}
\usetikzlibrary{automata, positioning}
\usepackage{mathtools}
\usepackage[rm={proportional,lining},sf={proportional,lining},tt={tabular,lining,monowidth}]{cfr-lm}
\renewenvironment{quotation}
               {\list{}{\listparindent=0pt
                        \itemindent
                        \listparindent
                        \leftmargin=16pt
                        \rightmargin=16pt
                        \topsep=6pt
                        \itemsep=0pt
                        \parsep=\medskipamount
                       }
                \item\relax}
               {\endlist}

\newcommand{\rel}{\ensuremath{\succsim}}

\newcommand{\U}{{\ensuremath{\mathscr{U}}}}
\newcommand{\Ux}{{\ensuremath{\mathscr{V}}}}
\newcommand{\sS}{{\ensuremath{\mathscr{S}}}}
\newcommand{\A}{{\ensuremath{\mathscr{A}}}}
\newcommand{\Ps}{{\ensuremath{\Pi}_F}}
\newcommand{\Uf}{{\ensuremath{\mathscr{U}_F}}}
\newcommand{\mR}{{\ensuremath{\mathbf{r}}}}
\newcommand{\mQ}{{\ensuremath{\mathbf{Q}}}}
\newcommand{\bP}{{\sbweight A1}}
\newcommand{\bPx}{{\sbweight A1}$'$}
\newcommand{\bA}{{\sbweight A2}}

\newcommand{\bC}{{\sbweight A3}}

\newcommand{\bUx}{{\sbweight A5}$'$}
\newcommand{\bS}{{\sbweight A4}}
\newcommand{\bU}{{\sbweight A5}}
\newcommand{\bComp}{{\sbweight A6}}
\newcommand{\bCont}{{\sbweight A7}}
\newcommand{\bCP}{{\sbweight A8}}
\newcommand{\bCPx}{{\sbweight A8}$'$}

\DeclarePairedDelimiter\abs{\lvert}{\rvert}%
\DeclarePairedDelimiter\norm{\lVert}{\rVert}%

\makeatletter
\let\oldabs\abs
\def\abs{\@ifstar{\oldabs}{\oldabs*}}
\let\oldnorm\norm
\def\norm{\@ifstar{\oldnorm}{\oldnorm*}}
\makeatother

\begin{document}

\title{An axiomatic approach to Markov decision processes}   
\author{Adam Jonsson}
\date{\today}
 \thanks{Tackar}

\institute{Adam Jonsson \at
                      Department of Engineering  Sciences and Mathematics
              \\Lule{\aa} University of Technology, 
              Lule{\aa}, Sweden\\ 
               \email{adam.jonsson@ltu.se}       
}

\date{Received: date / Accepted: date}

\maketitle

\begin{abstract}
This paper presents an axiomatic approach to finite Markov deci\-sion processes  where the discount rate is zero. One of the principal difficulties in the no discounting case is that, even if attention is restricted to stationary policies, a strong overtaking optimal policy need not exists. We provide prefe\-rence foundations for two criteria  that do admit optimal policies: $0$-discount optimality and average overtaking optimality. As a corollary of our results, we obtain  conditions on a decision maker's preferences which ensure that an optimal policy exists. These results have implications for  disciplines where  stochastic dynamic programming problems  arise, including   automatic control, dynamic games, and economic development.  
\end{abstract}

\section{Introduction}
This paper presents an axiomatic approach  to finite Markov decision processes (MPDs)  where the discount rate is zero.   MDPs  comprise a broad class of  stochastic dynamic deci\-sion problems and they have been studied extensively over the past several decades. To keep the discussion as elementary as possible, we will work within the framework of Blackwell's  \cite{Bla62} classic paper.  For extensions of this framework and discussion of its many uses, the reader is  referred to \cite{ABF93,HV98,RSV02} and the books \cite{FS02,Piu13,Put94}.

 In its simplest form, a  MDP has  the following ingredients: A state space $\sS$, an action space $\A$, a transition probability function $p_a(s' \vert s)$ on $\sS$ for each   $a \in \A$, and a real-valued function $r(s, a)$  on $\sS \times \A$. Here  $\sS$ represents possible states of a system (a manufacturing chain, a biological system,  a natural resource, etc.) and $\A$ represents choices available to an agent (the decision maker). Unless stated other\-wise, $\sS$ and $\A$ are    finite sets.  At discrete times $t=1, 2, 3, \ldots$, the agent observes the state and selects an element from $\A$. If the system is in   $s \in \sS$ and $a \in \A$ is chosen, then a reward of $r(s, a)$ is received and the system moves to $s'$ with probability   $p_a(s' \vert s)$. Rewards are discounted so that a reward of one unit at time $t$ has present value  $\beta^t$, where $0<\beta \leq 1$. The  problem is to choose a policy (i.e., a rule for selecting actions at all future times) that maximizes the expected net present value of all future rewards.

This problem is particularly   difficult when $\beta=1$.  To begin with, it is not  clear what it means to maximize net present value in this case. The difficulty is that the total  value of a policy is typically infinite if $\beta=1$.  There is a  natural sense in which a policy is maximal if it generates a sequence of cumulative  expected rewards that eventually  dominates that of any other policy. This leads to the intuitive notion of over\-taking optimality (formally defined in Section \ref{sect: example}).  It is well known, however, that an overtaking optimal  policy need not exist. A  less selective criterion is based on the expected long-run average reward of a policy. But this   criterion  does not differentiate  between  streams of expected rewards which might have very different appeal to the decision maker.

 Blackwell \cite{Bla62}  introduced the   \emph{$1$-optimality} criterion,  which evaluates streams of expected rewards  on the basis of their Abel means. He also established the existence of $1$-optimal policies that are \emph{stationary}, (i.e., for which  the action chosen at time $t$ depends only on the state of the system at time $t$).\footnote{More precisely, Blackwell  \cite{Bla62}  establishes existence  of optimal policies using   the criterion now known as \emph{Blackwell optimality}, which is slightly stronger than $1$-optimality.  He refers to $1$-optimality as \emph{near optimality}; other authors use the terms \emph{$0$-discount optimality}  and  \emph{bias optimality} \cite{Put94,Piu13}.} Subsequently, Veinott \cite{Vei66} introduced what is often referred to as the \emph{average overtaking} criterion,  where  Abel means are  substituted for {C}es\`{a}ro  means. The Blackwell--Veinott criteria are able to select between policies that the  average reward criterion 
does not distinguish. However, the literature has not adressed the following questions:

\begin{quotation}  
 \text{Q1}. Are the Blackwell--Veinott criteria the \emph{only} selective criteria which admit optimal policies in the no discounting case? 

 \text{Q2}.  How can these criteria  be described axiomatically?

 \text{Q3}. Under which assumptions on a decision maker's preferences do optimal policies exist?
\end{quotation}
Our main results are summarized in  Theo\-rem \ref{theorem1},  \ref{theorem2} and \ref{theorem3}. Theo\-rem \ref{theorem1}   shows that, subject to certain constraints,   \text{Q1} has an affirmative answer. Theorem  \ref{theorem2} and  \ref{theorem3}  provide two sets of axioms that characterize the average overtaking and $1$-optimality criterion on the reward streams generated by stationary policies. The second of these two results complements a theorem of Jonsson and Voorneveld  \cite{JV18} and uses the compensation principle  as a key axiom.  Finally,  we obtain a partial answer to    \text{Q3}    as a corollary of these results.  


\section{Preliminaries}\label{sect 2}
 Our finite MDP has state space $\sS$ and  action space $\A$. At  times $t=1, 2, 3, \ldots$,  the agent observes the state of  and chooses an element $a$  from $\A$. We assume that this choice   depends on the history of the system only through its present state. Thus, the action chosen at time $t$ is an element of $F$, the set of all functions from  $\sS$ to $\A$.  Each $f\in F$ has a corresponding  transition matrix, $\mQ(f)$, and reward vector,   $\mR(f)$. In the notation from the introduction, if the system is in $s \in \sS$ and $f$ is used, then a reward   of $\mR(f)_s=r(s, f(s))$ is received and the system moves to $s'$ with probability  $\mQ(f)_{s, s'}=p_{f(s)}(s' \vert s)$.  Rewards may be interpreted, for example, as    payouts of a single good received by an infinitely lived consumer, or as the   utilities  of   future  generations.  
 
A \emph{policy} is a sequence $(f_1, f_2, f_3, \dots)$ in $F$. Using  policy $\pi=(f_1, f_2, f_3, \dots)$  means that, for each $t=1, 2, 3, \ldots$,  $f_t(s)$ is selected from $\A$ if the system is  in state $s$.  A policy is \emph{stationary} if using it implies that  the action chosen at time $t$  depends on the state of the system at time $t$, but not on $t$ itself. Formally, a   stationary policy can be written $(f, f, f, \ldots)$ for some $f \in F$.\footnote{More general definitions of the concepts of a policy and stationary policy allow for randomized actions (see, e.g., \cite[p.~22]{Put94}). Our results for  non-randomized (or pure)  stationary policies generalize trivially to randomized stationary policies.}  We denote the set of all policies by  $\Pi$ and the set of all stationary policies by  $\Ps$. 
  
Given an initial  state $s\in \sS$, the sequence of expected rewards  that  $\pi \in \Pi$ generates  is denoted   $u(s, \pi)$. If $\pi=(f_1, f_2, f_3, \dots)$ and $u=(u_1, u_2, u_3, \ldots)=u(s, \pi)$, then  
  \begin{align}\label{eq: generated stream}
u_1&=[\mR(f_1)]_s,\notag\\
u_t&=[\mQ(f_1) \cdot \ldots   \cdot \mQ(f_{t-1})\cdot \mR(f_t)]_s, \;  t\geq 2.
\end{align}  
Let $\Uf$ be the set of sequences generated by stationary policies. That is, $u \in \Uf$ if and only if $u=u(s, \pi)$ for some $s\in \sS$ and $\pi \in \Ps$.   

 The agent needs to compare $u(s, \pi)$ and $u(s, \pi')$ for different $s \in \sS$ and $\pi, \pi' \in \Pi$. For convenience, we consider (incomplete) preferences on the set of all bounded sequences, which is denoted by $\U$. We  reserve the notation $\succsim$ for a preorder on $\U$ (i.e., a reflexive and transitive binary relation), where  $u \succsim v$ means  that $u$ is at least as good as $v$.   We  say that $\succsim$  \emph{compares} $u$ and $v$ if   either   $u \succsim v$ or $v \succsim u$, and we write $\neg u \succsim v$  to indicate that $u$ is not at least as good as $v$. As usual, $u \succ v$ denotes strict preference ($u \succsim v$, but $\neg v \succsim u$) and $u \sim v$ denotes indifference ($u \succsim v$ and $v \succsim u$).   

In this framework, preferences are thus defined over sequences of expected rewards. That is, it is assumed that preferences over random rewards can be reduced to preferences over expected rewards. The framework is therefore unable  to elucidate risk-averse prefe\-rences.  For risk measures and risk-sensitive control of Markov processes, see \cite{BR14,Rus10} and the references cited there. 


\section{A motivating example}\label{sect: example}
For background, we begin by  reviewing  how different ways of comparing reward streams  may fail or succeed to yield optimal policies. The comparisons  often involve  sums over a finite horizon.  
For  $u \in \U$ and $T \in \mathbb{N}$, we let 
\begin{align}\label{eq: sigma 0}
\sigma_T(u)=\sum_{t=1}^Tu_t,  
\quad  
\sigma(u)=(\sigma_1(u), \sigma_2(u), \sigma_3(u), \ldots).
\end{align} 

A policy $\pi^\ast \in \Pi$ is  \emph{overtaking optimal} if, for every $\pi \in \Pi$,
\begin{align}\label{eq: definition opt 0} 
& u(s, \pi^\ast)  \succsim_{\text{O}} u(s, \pi) \text{ for every }s\in \sS,
\intertext{where} \label{eq: definition OT} 
   & u  \succsim_{\text{O}}  v \; \Longleftrightarrow  \; \liminf_{T\to \infty} \sigma_T(u-v) \geq 0. 
\end{align} 
This criterion has the advantage of being plausible intuitively. It is also the strongest among the most commonly discussed criteria for undiscounted MDPs.  Its drawback is that an optimal policy need not exist  \cite{Bro65,Gal67}. The following  is a variation of an example from   Denardo and Miller \cite{DM68}. We return  to this example in  Section \ref{sect: characterizations}. 

\begin{example}\label{ex: 1}
Figure \ref{fig: 1} displays the transition graph of a deterministic MDP with $\A=\{a_1, a_2\}$ and $\sS=\{s_1, s_2, s_3\}$.   
If the system starts in state $s_1$ and  $a_1$ is chosen, then the system moves to $s_2$ and  a reward of $2$ is received; if $a_2$ is chosen,  the system moves to $s_3$  and  a reward of $c \in \mathbb{R}$ is received. Once the system reaches $s_2$ or $s_3$, it starts to alternate between these two states,   and it does not matter how the agent acts. A reward of $0$ is received when the system goes from $s_2$ to $s_3$, and a reward of $2$ is received when it goes from $s_3$ to $s_2$.  
\begin{figure}[h]
 \centering
\begin{tikzpicture}[font=\sffamily]
\label{fig: 1}
        \tikzset{node style/.style={state, minimum width=1cm, line width=0.3mm, fill=white!20!white}}
        \node[node style] at (0, 0) (s1) {{\large $s_1$}};
        \node[node style] at (6, 0) (s2) {{\large $s_2$}};
        \node[node style] at (3, -1.596) (s3) {{\large $s_3$}};

        \draw[every loop, auto=right, line width=0.25mm, >=latex, draw=black, fill=black]
        (s1) edge[bend right=20] node {$r(s_1, a_2)=c$} (s3)
        (s1) edge[bend left=20] node {$r(s_1, a_1)=2$} (s2)
        (s2) edge[bend right=20, auto=right] node {$r=0$} (s3)
        (s3) edge[bend right=20] node {$r=2$} (s2);
        
\end{tikzpicture}	
\caption{A deterministic MDP where no overtaking-optimal policy exists.}
\end{figure}

Suppose that the system starts in   $s_1$.  Let $u$ be the reward stream  that is generated if $a_1$ is chosen, and let $v$ be the stream that obtains if  $a_2$ is chosen.  Then 
\begin{align*} 
u=(2, 0, 2, 0, 2, \dots)  \quad  \text{ and } \quad v=(c, 2, 0, 2, 0, 2, \dots).
\end{align*} 
We have $\sigma_T(u-v)=2-c$ if $T$ is odd and $\sigma_T(u-v)=-c$ if $T$ is even. Hence,  if $0<c<2$, then $\neg u \succsim_{\text{O}}  v$ and $\neg v \succsim_{\text{O}}  u$. This means that there is no overtaking-optimal policy if $0<c<2$. 
\qed
\end{example}

Note that the MPD in Example  \ref{ex: 1} still does not admit an overtaking optimal policy if attention is restricted to  stationary policies. We remark that it is not only in deterministic MDPs that this limitation of   overtaking  optimality makes itself known. There are, indeed,  ergodic MDPs  where no overtaking-optimal policy exists within the class of stationary policies  \cite{NVA99}. 

 Let us also note that  optimal policies often do exist  if we adopt an alternative definition of overtaking optimality, according to which $\pi^\ast \in \Pi$  is optimal if    there is no $\pi \in \Pi$ such that  
\begin{align*} 
u(s, \pi) \succ_{\text{O}} u(s, \pi^\ast) \text{ for every }s\in \sS.
\end{align*}
(In Example \ref{ex: 1}, all policies are optimal in this sense if $0<c<2$.) This weaker form of overtaking optimality has been used frequently in   studies of optimal economic growth  \cite{Bro70a,BM73,BM07}. It is closely related to the notion of \emph{sporadic overtaking optimality}  studied in the operations research literature  \cite{Ste84,FPS17}. Here we have adopted the definition of overtaking optimality that this literature most frequently employs.  
 
Generalizing the definition   \eqref{eq: definition OT} to an arbitrary preorder $\succsim$,  let us say that $\pi^\ast \in \Pi$ is \emph{$\succsim$-optimal} or \emph{optimal with respect to $\succsim$}  if, for every $\pi \in \Pi$,
\begin{align}\label{eq: definition opt 1}
 u(s, \pi^\ast) \succsim  u(s, \pi) \text{ for every }s\in \sS.
\end{align}
The preorders associated with average reward optimality,  average overtaking optimality and $1$-optimality are defined as follows:

\begin{align}
  \label{eq: definition AV}
& \textbf{ (average reward) }  && u \succsim_{\text{AR}} v \Longleftrightarrow   \liminf_{T\to \infty}\frac{1}{T}\sigma_T(u-v) \geq 0 \\
 \label{eq: definition AOT}
 &\textbf{ (average overtaking) } && u \succsim_{\text{AO}} v \Longleftrightarrow   \liminf_{T\to \infty} \frac{1}{T}\sum_{t=1}^T\sigma_t(u-v) \geq 0\\
  \label{eq: definition 0-disc}
 &\textbf{ ($1$-optimality) }
  && u \succsim_1 v \Longleftrightarrow   \liminf_{\delta \to 1^{-}}  \sum_{t=1}^\infty \delta^{t}\cdot (u_t-v_t)\geq 0.
\end{align}
The average reward criterion is the most studied criterion for undiscounted MDPs. The standard criticism against this criterion concerns the fact that improvements in any finite number  of time periods are ignored. In Example  \ref{ex: 1},  for instance, it is average reward-optimal to choose $a_1$ in state $s_1$ even if the value of  $c$ is very large.

If $u$ and $v$ are the  streams in Example \ref{ex: 1},  the  {C}es\`{a}ro sum of $\sum_{t=1}^\infty(u_t-v_t)$  is $1-c$. Hence, it is average overtaking-optimal to choose $a_1$  if and only if $c\leq 1$.  It is well known that  average overtaking optimality  is equivalent to  $1$-optimality  in finite MDPs \cite{Lip69}.  In general, any average overtaking-optimal policy is $1$-optimal, but a $1$-optimal policy need not be  average overtaking optimal (see, e.g., \cite{BFZ14}).

To sum up, while the average reward criterion is unselective, the overtaking criterion is overselective.  One way to formulate the first question (\text{Q1}) from the introduction is to ask 
if the average overtaking criterion is the least selective criterion that admits optimal policies.  To state this question in a precise way, we will formulate a set of conditions which we can plausibly require of a selective criterion. 


   \section{Axioms}\label{sect: axioms}
This section provides five conditions (called axioms) on preorders that are  known from the literature. The five conditions are satisfied by the preorders associated with the overtaking criterion, the average overtaking criterion and the $1$-optimality criterion (see \cite[p.~28]{JV18}). They may be viewed as conditions that can be plausibly required of a selective criterion.  
 
The first axiom, {\bP}, is a standard monotonicity requirement. It asserts that preferences are positively sensitive to improvements in each time period.  
Preorders that  meet this requirement avoid the standard criticism of the average reward criterion.
 \begin{quotation}
 \bP.  For  all $u,v \in \U$, if $u_t\geq v_t$ for all $t$  and  $u_t>v_t$ for some $t$, then $u \succ v$.
 \end{quotation}
This axiom says, in particular, that the agent  prefers a certain reward of $2$ units to a certain reward of $1$ unit.   In the present framework,  it also says that the agent disprefers a certain reward of $2$ units to a lottery that pays a reward of $1$ or $4$ units with equal probabilities. As indicated in Section \ref{sect 2}, such assumptions are inappropriate for risk-averse agents.

 The second axiom, {\bA}, formalizes the assumption that    a reward of one unit at time $t>1$ is worth the same as   a reward of one unit at  $t=1$  (i.e., that $\beta=1$). In the case when   rewards represent utilities (or consumption) of future generations,   {\bA} is  the axiom of \emph{anonymity}, which ensures the equal treatment of generations.   \begin{quotation}
\bA.  For all $u,v \in \U$, if $u$ can be obtained from $v$ by interchanging two entries of $v$, then  $u \sim v$.  
 \end{quotation}

 The next axiom is a relaxation of the consistency requirement used in Brock's  \cite{Bro70b} characterization of the overtaking criterion. For $n \geq 1$ and $u \in \U$, let  $u_{[n]}$  be   the stream obtained from $u$ by replacing $u_t$ with $0$ for all $t>n$. Our third axiom can then be stated as follows. 
  \begin{quotation}
\bC.  For  all $u,v \in \U$, if there exists $N > 1$ such that $u_{[n]} \succ v_{[n]}$ for all $n\geq N$, then  $u \succsim v$.
\smallskip 
 \end{quotation}
 That the average reward criterion   satisfies {\bC} is a trivial consequence of the fact that $u_{[n]} \sim_{\text{AR}} v_{[n]}$ for all $u, v \in \U$  and every $n\geq 1$. The preorders in \eqref{eq: definition OT}, \eqref{eq: definition AOT} and \eqref{eq: definition 0-disc}  have the stronger property that   $u$ is at least as good as $v$ if  $u_{[n]}$ is merely at least as good as $v_{[n]}$ for all sufficiently large $n$; this property fails for the average reward criterion.   

The fourth   axiom asserts that for reward streams $u, v \in \U$, if both streams are postponed one period and an arbitrary reward  of $c \in \mathbb{R}$ is assigned to the first period, then the resulting streams, $(c, u)=(c, u_1, u_2, u_3, \ldots)$ and $(c, v)=(c, v_1, v_2, v_3, \ldots)$,  should be ranked in the same way as $u$ and $v$.
 \begin{quotation}
\bS. For  all $u,v \in \U$ and  $c \in \mathbb{R}$, $(c, u)  \succsim (c, v)$  if and only if $u  \succsim v$. 
\end{quotation}
This axiom was proposed as a fundamental  condition by Koopmans \cite{Koo60} in his pioneering work on  intertemporal choice. It is usually referred to as  \emph{stationarity} \cite{AdB10,BRW08} or \emph{independent future} \cite{FM03,Mit18}.

Our last  axiom is an adaptation of the standard assumption of interpersonal comparability from social choice theory (see, e.g., \cite{DG77}). In the intertemporal setting, it asserts that preferences are invariant to changes in the origins of the utility   indices used in different  periods. This condition has been referred to as   \emph{zero  independence} \cite{Mou88} and   
 \emph{translation scale invariance}  \cite{AdB10}.
 \begin{quotation}
\bU.  For all $u,v, \alpha \in \U$, if  $u \succsim v$, then $u+\alpha \succsim v+\alpha$. 
\smallskip 
\end{quotation}
Note that a preorder $\succsim$ 
which satisfies {\bU} has the property that if $u, v, u', v' \in \U$ are such that $u-v=u'-v'$, then $u \succsim v$ if and only if $u' \succsim v'$.  (The converse is also true.) This fact will be used repeatedly below.


\section{A rigidity result}
If we view the axioms from the previous section   as conditions which we expect a selective criterion to satisfy, then the first question from the introduction can  be stated as follows:  
 If  $\succsim$  satisfies {\bP}--{\bU}, is every $\succsim$-optimal policy  average overtaking-optimal (and hence $1$-optimal)?\footnote{An alternative way to state   \text{Q1}  would be to ask if  $\succsim_{\text{AO}}$ is the least restrictive extension of $\succsim_{\text{O}}$ that admits optimal policies. This question has a trivial answer,  however, because $\succsim_{\text{AO}}$ is not, strictly speaking, even an extension of $\succsim_{\text{O}}$: if $u \succsim_{\text{O}} v$, then  $u \succsim_{\text{AO}} v$, but there are $u, v \in \U$ with  $u \succ_{\text{O}} v$ and $u \sim_{\text{AO}} v$ (see \cite[p.~28]{JV18}).} 
Theorem \ref{theorem1}
 shows that this question has an   affirmative answer if  attention is restricted to stationary policies.   This restriction does not trivialize any of the questions (Q1--Q3) from the introduction. In fact, replacing $\Pi$ with  $\Ps$ in the preceding discussion would not affect what has been said so far in an essential way.

\begin{theorem}\label{theorem1}
Suppose that $\succsim$  satisfies {\bP}--{\bU}. If a policy is $\succsim$-optimal within the class of stationary policies, then it is average overtaking-optimal  within the class of stationary policies.
\end{theorem}

\begin{proof} 
 The proof   exploits the fact that under certain conditions on $u \in \U$, if a preorder $\succsim$ satisfies {\bP}--{\bU}, then 
\begin{align}\label{eq: compensation implication1}
u \succsim  (0, u) \text{ implies } \bar{u} \geq 0, 
\end{align} 
where 
\begin{align}\label{def: average}
   \bar{u}\equiv \lim_{n\to \infty}\frac{1}{n}\sum_{t=1}^nu_t
\end{align}
is the average of $u$. The usefulness of \eqref{eq: compensation implication1} is explained by the fact that if $\succsim$ satisfies {\bU} and $u, v \in \U$ are such that $\sigma \equiv \sigma(u-v)$ is bounded, then 
\begin{align}\label{eq: bU and equivalence}
\text{ $u \succsim v$ if and only $\sigma \succsim (0, \sigma)$.}
\end{align}
This is because $u-v=\sigma-(0, \sigma)$. Applying \eqref{eq: compensation implication1} with $\sigma$ in the role of $u$, we see that $u \succsim v$ implies $\bar{\sigma}\geq 0$. Since $\bar{\sigma}$ is the  {C}es\`{a}ro sum of $\sum_{t=1}^\infty(u_t-v_i)$, this means that $u \succsim v$ implies $u \succsim_{\text{AO}} v$.  

The conditions on $u \in \U$  which ensure   \eqref{eq: compensation implication1}  
are that (i)  the limit \eqref{def: average}  
exists and (ii) for every $\varepsilon>0$ there exists an $N$ such that the average of any $n\geq N$ consecutive coordinates of $u$ differs from $\bar{u}$ by at most $\varepsilon$---that is,   
\begin{align*}
   \abs{ \frac{1}{n}\sum_{t=t_0}^{t_0+n}u_t-\bar{u}} <\varepsilon \text{ for every   }t_0 \in \mathbb{N}. 
\end{align*}
 We say that $u \in \U$ is \emph{regular} if the two conditions are met.

\begin{lemma}{\cite[Proposition 1]{JV18}}
\label{lemma1}
Suppose that  $\succsim$ satisfies {\bP}--{\bU}. If $u \in \U$ is regular and $c \in \mathbb{R}$, then 
\begin{align*} 
(c, u) \rel u \text{ implies }c\geq \bar{u} \intertext{and} 
u \rel (c, u) \text{ implies }c\leq \bar{u}.
\end{align*} 
 
\end{lemma}

Now, for every $\pi \in \Ps$,  $u(s, \pi)$ is regular  for each  $s \in \sS$. This follows from the well known fact that  the  reward stream generated by a stationary policy can be written as the sum of a periodic sequence and a summable sequence. (The stream generated by $(f, f, f, \ldots)$ is defined by powers of  $\mQ(f)$ acting on $\mR(f)$---see \eqref{eq: generated stream}. By the Perron-Frobenius theorem for non-negative matrices, the sequence $\mQ(f) \cdot \mR(f), \mQ(f)^2\cdot \mR(f), \mQ(f)^3\cdot \mR(f), \ldots$ approaches a periodic orbit at exponential rate.)    To apply the arguments preceding Lemma \ref{lemma1}, we need to know that $\sigma(u-v)$ is bounded and regular if $u$ and $v$ are generated by   stationary policies.  We have the following result.

\begin{lemma}\label{lemma2}
Suppose that $u$ and $v$ are generated by stationary policies, and let $\sigma \equiv \sigma(u-v)$ be defined as in \eqref{eq: sigma 0}.  If $\bar{u}=\bar{v}$, then $\sigma \in \U$ is regular. 
\end{lemma}
\begin{proof}
Write  
\begin{align}\label{eq: decomposition}
u=x^{(u)}+y^{(u)},  \quad v=x^{(v)}+y^{(v)},
\end{align}
where $x^{(u)}$ and $x^{(v)}$ are periodic and where $y^{(u)}$ and $y^{(v)}$ are summable. Let $p$ be the product of the periods of $x^{(u)}$  and $x^{(v)}$. Then $\bar{u}=\bar{x}^{(u)}=\sigma_p(x^{(u)})/p$ and $\bar{v}=\bar{x}^{(v)}=\sigma_p(x^{(v)})/p$. So, if $\bar{u}=\bar{v}$, then $\sigma_p(x^{(u)}-x^{(v)})=0$. This means that $\sigma(x^{(u)}-x^{(v)})$ is periodic. The sequence  $\sigma(y^{(u)}-y^{(v)})$ is convergent by our choice of $y^{(u)}$ and $y^{(v)}$. Hence,  $\sigma=\sigma(u-v)$ is the sum of a periodic   sequence and a convergent sequence. This means that $\sigma  \in \U$ is  regular.   \qed 
\end{proof}

To complete the proof of  Theorem  \ref{theorem1}, let $\succsim$ be a preorder that satisfies {\bP}--{\bU}, and suppose that $\pi^\ast$ is $\succsim$-optimal within $\Ps$. Let $u=u(s,\pi^\ast)$ and $v=u(s,\pi)$, where $\pi \in \Ps$ and $s \in \sS$ are  arbitrary, and let  $\sigma\equiv \sigma(u-v)$ be defined as in \eqref{eq: sigma 0}.  Since  $\pi^\ast$ is $\succsim$-optimal within $\Ps$,   $u \succsim v$. We need to show that  $u \succsim_{\text{AO}} v$. If $\bar{u}=\bar{v}$, then this follows from Lemma \ref{lemma1} and \ref{lemma2}   and the remarks preceding Lemma \ref{lemma1}. It remains  to show that  $u \succsim_{\text{AO}} v$  if  $\bar{u} \neq \bar{v}$. It is enough to show that $\bar{u} > \bar{v}$, since this clearly implies $u \succ_{AO} v$.  Given any preorder $\succsim'$ that satisfies {\bP}--{\bU},   if $x \in \U$ and $y\in\U$ are such that $\bar{x} > \bar{y}$, then $x \succ' y$ (see \cite{BM07} or \cite{JV15}).  Thus, if  $\bar{u} \neq \bar{v}$, then  we must  have  $\bar{u} > \bar{v}$. (If it were the case that $\bar{v} > \bar{u}$, then we would have $v \succ u$, which contradicts the assumption that $u \succsim v$.)  We can therefore conclude that $u \succ_{\text{AO}} v$, and the proof of Theorem  \ref{theorem1} is thereby complete.      \qed 
 \end{proof}
 
 
 \section{Characterizations}\label{sect: characterizations}
One goal of this paper is to provide a preference foundation for finite MDPs.   In the case of a positive discount rate, the well known preference foundation of Koopmans \cite{Koo60,Koo72} is easily adapted to the present setting. The literature provides characterizations of two   criteria for the no discounting case: the over\-taking criterion  \cite{AT04,BM07,Bro70b} and the average reward criterion   \cite{KSW14,Mar98,KS18,Piv22}. The overtaking criterion is characterized by axioms that are similar to those   in Section \ref{sect: axioms}. The characterizations of the average reward criterion, which does not satisfy {\bP},  involve further conditions of permutability and numeric representability.  These conditions are well known to be incompatible with {\bP} in the no discounting case \cite{BM03,FM03}.

In this section, we axiomatize the preorders associated with the average overtaking criterion and the $1$-optimality criterion.  As in the previous section, we restrict attention to stationary policies.  
 
  \subsection{First characterization}\label{subsect: 1st characterization}
The axioms from Section \ref{sect: axioms}  do not characterize $\succsim_{\text{AO}}$. Indeed, the preorder associated with the overtaking criterion satisfies  {\bP}--{\bU},  and $\succsim_{\text{O}}$ does  not agree with $\succsim_{\text{AO}}$ on $\Uf$ (the   streams generated by stationary policies). As illust\-rated in Example  \ref{ex: 1}, for $\succsim_{\text{AO}}$-optimality to imply $\succsim$-optimality, it is necessary that $\succsim$  compares at least some pairs of streams   that $\succsim_{\text{O}}$ does not compare.  

Insisting that all pairs $u, v \in \U$  be comparable has unwanted consequences. In fact, it is not possible to give an explicit definition of a preorder, satisfying  {\bP} and {\bA}, that compares all pairs of sequences of $0$s and $1$s \cite{Lau10}.   On the other hand, $\succsim_{\text{AO}}$   compares each pair $u, v \in \Uf$ and coincides  with $\succsim_{1}$  
on this domain. Thus, the following condition is compatible with {\bP}--{\bU}:
 \begin{quotation}
{\bComp}. For all $u, v \in \Uf$,  $\succsim$  compares $u$ and $v$.
\smallskip 
\end{quotation}
If  $\succsim$ satisfies  {\bP}--{\bComp} and $u, v \in \Uf$,  then $u \succ v$ if and only if  $u \succ_{\text{AO}} v$. To conclude that the symmetric parts of $\succsim$ and $\succsim_{\text{AO}}$ agree,   further assumptions are needed. A sufficient condition asserts that, for all $u, v \in \U$, if $(\varepsilon+u_1, u_2, u_3, \ldots) \succsim v$ for every $\varepsilon>0$, then $u \succsim v$. This condition can be formalized by defining a metric on $\U$ and demanding that $\{v \in \U\colon u \succsim v\}$ be a closed subset of $\U$ for every $u \in \U$.  Almost any metric from the literature  will do (e.g., \cite[p.~5]{BM08}). For  example, let   $d(u, v)=\min\{1, \sum_{i=1}^\infty \vert u_i-v_i \vert \}$. The  continuity requirement can then be stated as follows.

 \begin{quotation}
{\bCont}. For every $u \in \U$, $\{v \in \U\colon u \succsim v\}$ is a closed subset of $\U$.
\smallskip 
\end{quotation}
 
\begin{theorem}\label{theorem2}
If $\succsim$  satisfies {\bP}--{\bCont}, then  $\succsim$ and $\succsim_{\text{AO}}$  coincide on $\Uf$. 
\end{theorem} 
\begin{proof}
Let $\succsim$  satisfy {\bP}--{\bCont}, and let $u, v \in \Uf$. We know that $u \succsim_{\text{AO}} v$ if $u \succsim v$  (Theorem \ref{theorem1}). So it is enough  to   show  that $u \succsim_{\text{AO}} v$ implies $u \succsim v$.

If $u \succ_{\text{AO}} v$, then either   (i) $\bar{u}>\bar{v}$ or (ii) $\bar{u}=\bar{v}$ and $\bar{\sigma}>0$,  where $\sigma=\sigma(u-v)$. In case (i), we get $u \succ  v$ as a consequence of the fact that $\succsim$  satisfies {\bP}--{\bU}. In case (ii), $\neg (0, \sigma) \succsim \sigma$  by Lemma \ref{lemma1}, so $\neg v \succsim u$ by {\bU}. By {\bComp}, $u \succ v$. Conclude that $u \succ_{\text{AO}} v$ implies $u \succ v$. 

 Now suppose that $u \sim_{\text{AO}} v$.  Let  $u^{(\varepsilon)}=(\varepsilon+u_1, u_2, u_3, \ldots)$.  Then $u^{(\varepsilon)} \succ_{\text{AO}} v$ for every  $\varepsilon>0$,  so (by the above conclusion) $u^{(\varepsilon)} \succ v$ for every  $\varepsilon>0$. By {\bCont}, $u \succsim v$. The same argument shows that $v \succsim u$.  \qed
 
\end{proof}
  
  \subsection{Second characterization}\label{subsect: 2nd characterization}
 Axioms {\bComp} and {\bCont} were motivated by necessity rather than some normative or economic reason. In our second characterization,  these axioms are replaced by the \emph{compensation principle}.   

As an illustration of this principle, imagine that the decision maker is faced with   two options. The first option yields some sequence of expected rewards $u\in \U$. 
The second option is to obtain a one-period postponement of $u$ and a compensation of $c \in \mathbb{R}$  in the first period. Which value of $c$ should make the agent indifferent? 

In some cases, this value will be zero. This is the case if $u$ has at most  finitely nonzero entries---then $(0, u)$ and  $u$ are equally good by {\bA}. However,  the agent will not always be indifferent if  $c=0$. For instance, if $u=(r, r, r, \ldots)$ is constant and $c$ is less than $r$, then $(c, u)$ is worse than $u$ by {\bP}. The  compensation principle says that 
$u$ and $(c, u)$ are equally good if $c=\bar{u}$  (compare Lemma \ref{lemma1}). Its precise statement  is  as follows: 
 \begin{quotation}
{\bCP}.   For every  $u\in \U$, if $\bar{u}$  is well defined, then $(\bar{u},u) \sim u$.
\smallskip 
\end{quotation}
 
For a case of the two options  described above, consider again the system in Figure \ref{fig: 1}, and suppose that the system starts in $s_1$. The agent then has two options. If $a_1$ is chosen, then   $u=(2, 0,2,0,2,\dots)$ obtains. If $a_2$ is chosen, then $(c, u)$ is  obtained. Thus, the two feasible alternatives are  $u$ and $v=(c, u)$. Since $\bar{u}=1$, {\bCP} says that $u$ and $v$ are equally good if $c=1$.  

Example \ref{ex: 1} illustrates  the fact that $\succsim_{\text{O}}$  violates  {\bCP}. It is easy to check that  $\succsim_{\text{AO}}$    satisfies {\bCP}, and the same is true of $\succsim_1$  \cite{JV18}. To see that the average reward criterion also satisfies {\bCP}, note that  if $d=(c, u)-u$, then we have  $\sigma_T(d)=c-u_T$ and therefore $\liminf_{T\to \infty }\frac{1}{T}\sigma_T(d)=\liminf_{T\to \infty }\frac{1}{T}\sigma_T(-d)=0$.  It follows that  $(c, u) \sim_{\text{AR}} u$ for \emph{every} $c \in \mathbb{R}$ and $u\in \U$.

Like \eqref{eq: compensation implication1}, the usefulness of  {\bCP} stems from the fact that if $\succsim$ satisfies {\bU} and $u, v \in \U$ are such that $\sigma \equiv \sigma(u-v)$ is bounded, then $u \succsim v$ if and only $\sigma \succsim (0, \sigma)$. Thus, if $\succsim$ satisfies {\bP}, {\bU} and {\bCP}, then $u \succsim v$ if and only $\bar{\sigma} \geq 0$. In \cite{JV18}, this observation is used to  characterize $\succsim _1$ on the set of streams that are summable or eventually periodic. Theorem  \ref{theorem3}  extends this result to  streams that can be decomposed according to \eqref{eq: decomposition}.  

\begin{theorem}\label{theorem3}
If $\succsim$ satisfies  {\bP}, {\bU} and {\bCP}, then   $\succsim$ and $\succsim_{\text{AO}}$  coincide on $\Uf$.   \end{theorem} 
\begin{proof}

Let $\succsim$ be a preorder that satisfies  {\bP}, {\bU} and {\bCP}. For $u, v \in \Uf$, let $\sigma=\sigma(u-v)$. Suppose  that $\bar{u}=\bar{v}$. Then $\sigma \in \U$ is regular (Lemma \ref{lemma2}), which means that $\bar{\sigma}$ is well defined.    By {\bP} and {\bCP}, $\sigma \succsim (0, \sigma)$ if and only if $\bar{\sigma} \geq 0$. By  {\bU},  $u \succsim v$ if and only if $\sigma \succsim (0, \sigma)$. Hence,  $u \succsim v$ if and only if $\bar{\sigma} \geq 0$. Since $\bar{\sigma}$ is the {C}es\`{a}ro sum of $\sum_{t=1}^\infty(u_t-v_i)$, we see that  $u \succsim v$ if and only if $u \succsim_{\text{AO}} v$.

Now suppose (without loss of generality) that $\bar{u}>\bar{v}$. Then  $u \succ_{\text{AO}} v$. We show that $u \succ v$.   For $T>1$, define $z \in  \U$ by setting $z_t = u_t$ for $t \leq  T$ and $z_t = u_t-c$ for $t > T$. Then $z$ is the sum of periodic sequence and a summable sequence, and $u \succ z$  by {\bP}.  Since  $\bar{u}>\bar{v}$, we can choose $T$ so that $\sigma_t(u-z)\geq 0$ for all $t \geq T$. Since $\bar{z}=\bar{v}$, the preceding argument gives that $z \succsim v$, so $u \succ v$ by transitivity.  \qed
\end{proof}
 
We can obtain a characterization of average overtaking optimality in general discrete time MDPs by generalizing {\bCP}. This result, which concerns optimality within the class of all policies, is provided in the appendix. There we also verify  that the axioms in Theorem \ref{theorem3}  are logically independent.  
 
Theorem \ref{theorem2} and \ref{theorem3} provide two axioms sets that characterize $\succsim_{\text{AO}}$  on $\Uf$.  As a corollary of these results, we obtain a partial answer to the third question (Q3) from the introduction: If  $\succsim$ satisfies the axioms in any one of these axiom sets, then a policy  is   $\succsim$-optimal within $\Ps$ if and only if it is  $\succsim_{\text{AO}}$-optimal   within $\Ps$. In particular, a $\succsim$-optimal policy  exists  within $\Ps$. 

\begin{acknowledgements}
Thanks to Henrik Hult for a discussion on the topic and to Eugene Feinberg for very helpful feedback on an earlier version of the paper. This  version has greatly benefitted  from comments by two anonymous referees. 
\end{acknowledgements}

 \section*{Appendices}
Appendix A contains a characterization result on average overtaking optimality within the set of all policies. Appendix B establishes that the axioms used in Theorem  \ref{theorem1} and  Theorem \ref{theorem3} are logically independent.  
 
  \section*{A Average overtaking optimality within the set of all policies}
 Theorem \ref{theorem4}  below  provides a characterization of average overtaking optimality   in general discrete time MDPs.  In particular, we  make  no assumptions on the state and action spaces.

 To allow for unbounded reward functions, let us substitute $\U$, the set of bounded sequences, for $\Ux=\mathbb{R}^\mathbb{N}$, the set of all  real sequences. The reward stream generated by a non-stationary policy need not be regular, so  its average may be undefined. For $u \in \Ux$, we let 
\begin{align}\label{def: limits}
 \bar{u}_\ast=\liminf_{n \to \infty}\frac{1}{n}\sum_{t=1}^nu_t, 
   \quad \bar{u}^\ast=\limsup_{n \to \infty}\frac{1}{n}\sum_{t=1}^nu_t.
\end{align}

Our characterization result for discrete time MDPs uses the following three properties of the average overtaking criterion:     
    \begin{quotation}
 \bPx.  For  all $u,v \in \Ux$, if $u_t\geq v_t$ for all $t$  and  $u_t>v_t$ for some $t$, then $u \succ v$.
   \smallskip 
   
\bUx.  For all $u,v, \alpha \in \Ux$, if  $u \succsim v$, then $u+\alpha \succsim v+\alpha$. 

\bCPx.  For every $u \in \Ux$,  if $\bar{u}^\ast$ is finite, then $(\bar{u}^\ast, u) \succsim u$.  If $\bar{u}^\ast =+\infty$, then $u \succsim (0, u)$. 
\smallskip 
\end{quotation}
That $\succsim_{\text{AO}}$ satisfies {\bCPx} is easy to see once we observe that for $u \in \Ux, c \in \mathbb{R}$, if $v=(c, u)$,  then we have $\sigma_t(u-v)=u_t-c$ and $\sigma_t(v-u)=c-u_t$. Hence, 
\begin{align*}
\liminf_{n \to \infty}\frac{1}{n}\sum_{t=1}^n\sigma_t(v-u)= c-\bar{u}^\ast. 
\end{align*}
It follows that $(\bar{u}^\ast,  u)  \succsim_{\text{AO}}  u$ if  $\bar{u}^\ast$ is finite, and that $u \succsim (0, u)$ if $\bar{u}^\ast =+\infty$. 

Let us also note that every preorder that satisfies {\bCPx} and {\bUx}  also has the property that for all $u \in \Ux$,  
\begin{align}\label{eq: alternative CPx}
 \text{ if $\bar{u}_\ast$ is finite, then } u \succsim  (\bar{u}_\ast, u).
\end{align}
To see this, let $x=-u$. Then $\bar{x}^\ast=-\bar{u}_\ast$, so {\bCPx} implies $(-\bar{u}_\ast, -u) \succsim -u$. By {\bUx}  (adding $\alpha=u+ (\bar{u}_\ast, u)$), this means that $u \succsim (\bar{u}_\ast, u)$. In particular,  
\begin{align}\label{eq CPs}
 u  \succsim_{\text{AO}}  (\bar{u}_\ast,  u) \quad  \text{ and }  \quad (\bar{u}^\ast,  u)  \succsim_{\text{AO}}  u.
\end{align}
 
\begin{theorem}\label{theorem4} Let $\succsim$ be a preorder on $\Ux$ that satisfies  {\bPx}, {\bUx} and {\bCPx}. If a policy is $\succsim_{\text{AO}}$-optimal, then it is also $\succsim$-optimal. 
\end{theorem}
Note that Theorem \ref{theorem4} concerns the implication from $\succsim_{\text{AO}}$-optimality to $\succsim$-optimality whereas  Theorem \ref{theorem1} concerns the reverse implication.  As indicated above, Theorem \ref{theorem1}  does not hold in  non-finite MDPs. For example, the $1$-optimality criterion satisfies  {\bPx}, {\bUx} and {\bCPx}, and a $1$-optimal policy need not be average overtaking optimal. Whether or not Theorem \ref{theorem1}  holds in finite MDPs, without  restricting to stationary policies, is a question that we have not been able to answer. 
\begin{proof}   
Let $\succsim$ be a preorder on $\Ux$ that satisfies  {\bPx}, {\bUx} and {\bCPx}. Let $u \in \Ux$ be the stream of expected rewards generated by a  $\succsim_{\text{AO}}$-optimal policy, given some initial state, and let $v \in \Ux$ be generated by some other policy for the same initial state.  We need to show that $u  \succsim v$. 

Let $\sigma_n=\sigma_n(u-v), n \geq 1$, and let $\sigma=(\sigma_1, \sigma_2, \sigma_3 \ldots) \in \Ux$. That $u$ is generated by a  $\succsim_{\text{AO}}$-optimal policy means that  $u \succsim_{\text{AO}} v$.  By the definition of   $\succsim_{\text{AO}}$,  this implies that $\bar{\sigma}_\ast \geq 0$. Suppose first that $\bar{\sigma}_\ast < +\infty$. Since $ \succsim$ satisfies  {\bUx} and   {\bCPx}, we then have  $\sigma \succsim (\bar{\sigma}_\ast, \sigma)$ (see \eqref{eq: alternative CPx}).  By {\bPx} and transitivity, we thus have $\sigma \succsim (0, \sigma)$. By {\bUx} and the fact that  $u-v=\sigma-(0, \sigma)$,   this entails $u  \succsim v$.  If $\bar{\sigma}_\ast$ equals $+\infty$, then so does  $\bar{\sigma}^\ast$.  We then have $\sigma \succsim (0, \sigma)$ by {\bCPx} and hence $u  \succsim v$ by {\bUx}.  Conclude that $u  \succsim v$. Since $v$ was generated by an arbitrary policy, this shows that any  $\succsim_{\text{AO}}$-optimal policy is $ \succsim$-optimal.  \qed 
 
\end{proof}

\section*{B Logical independence}
\subsection*{Independence of the axioms in Theorem  \ref{theorem1}}
The following  binary relations, defined for all $u, v \in \U$, fail to satisfy precisely one of the axioms used in Theorem  \ref{theorem1}:  
\begin{align*}
  u \succsim_{\neg \text{\bP}} v \Longleftrightarrow  & u, v \in \U \text{ (all streams are equivalent) }  \\
  u \succsim_{\neg \text{\bA}} v \Longleftrightarrow  & \sum_{t=1}^\infty 2^{-t}(u_t-v_t)\\
       u \succsim_{\neg \text{\bC}} v \Longleftrightarrow & \exists T_0 \in \mathbb{N} \text{ s.t. } u_t \geq v_t \text{ for all $t > T_0$ and }   \sigma_{T_0}(u-v) \geq 0 \\[2mm]
              u \succsim_{\neg \text{\bS}}  v \Longleftrightarrow &  \liminf_{T\to \infty}   \sigma_{2\cdot T}(u-v) \geq 0  \, \text{ (cf.  \cite[p.~786]{FM03})}\\[2mm]
     u \succsim_{\neg \text{\bU}} v \Longleftrightarrow  &  \liminf_{T\to \infty}    \sigma_T(u^3-v^3) \geq 0.
\end{align*}
We omit the proofs for the first three  of these five preorders. The fourth  clearly satisfies {\bP}, {\bA} and {\bU}. To verify {\bC},  note that for any $u, v \in \U$ and $n \in \mathbb{N}$, we have $u_{[n]} \succ_{\neg \text{\bS}}  v_{[n]}$ if and only if $\sigma_n(u-v) > 0$. (Recall that  $u_{[n]}$  is the stream obtained from $u$ by replacing $u_t$ with $0$ for   $t>n$.)  Thus, if  $u_{[n]} \succ_{\neg \text{\bS}}  v_{[n]}$ for all sufficiently large $n$,   then $\sigma_n(u-v) > 0$ for all sufficiently large $n$.   In particular, $\sigma_{2\cdot n}(u-v) >0$ for all sufficiently large $n$, which means that $u \succsim_{\neg \text{\bS}}  v$.  
 Conclude that {\bC} holds. 

 To see that $\succsim_{\neg \text{\bS}}$ violates {\bS}, let $u$ and $v$ be the  streams in Example \ref{ex: 1} with $c=1$, so that $u=(2, 1, 2, \ldots)$ and $v=(1, 2, 1, 2, \ldots)$. By  the definition of $\succsim_{\neg \text{\bS}}$, we  have $u \sim_{\neg \text{\bS}}  v$ and $(2, u)  \succ_{\neg \text{\bS}}  (2, v)$. This shows that $\succsim_{\neg \text{\bS}}$ fails to satisfy {\bS}.

It is  straightforward to check that $\succsim_{\neg \text{\bU}}$ satisfies {\bP}--{\bS}. To show that  {\bU} fails, let  $u = (3, 0, 0, 0, \ldots), v = (2, 2, 0, 0, 0, \ldots)$,  
$\alpha = (0, 1, 0, 0, 0, \ldots)$. Define $x= u+\alpha$, $y=v+\alpha$. Then  $x=(3, 1, 0, 0, 0, \ldots)$ and $y=(2, 3, 0, 0, \ldots)$.  So, for $T \geq 2$, we have $\sigma_T(u^3-v^3)=11$ and  $\sigma_T(x^3-y^3)=-7$. Hence, $u \succsim_{\neg \text{\bU}} v$, but $u+\alpha \neg  \succsim_{\neg \text{\bU}} v+\alpha$. This shows that {\bU} fails. 
 
The five preorders    ($\succsim_{\neg \text{\bP}}$ to $\succsim_{\neg \text{\bU}}$) establish logical independence of  {\bP}-- {\bU}. The first, second, fourth and fifth  preorder   show that  Theorem \ref{theorem1} fails if we drop any one of    {\bP},  {\bA},  {\bS} and {\bU}.  We have been unable to find a preorder which shows that the theorem fails if {\bC} is dropped.  
 
\subsection*{Independence of  the axioms in Theorem  \ref{theorem3} }
We show that the axioms in Theorem  \ref{theorem3}  are logically independent by providing three preorders, each violating precisely one of the three axioms. The overtaking criterion    
satisfies {\bP} and {\bU}, but not {\bCP} (see Example \ref{ex: 1}). 
The preorder $\succsim_{\neg \text{\bP}}$  satisfies {\bU} and {\bCP}, but not {\bP}.  (The average reward criterion provides another example.)  It   remains to find a preorder that satisfies {\bP} and {\bCP}, but violates {\bU}. 

For $u, v \in \U$,   let us we say that $u$ \emph{dominates} if  $u_t\geq v_t$ for all $t \in \mathbb{N}$, and that $u$ is a \emph{finite permutation} of $v$ if $u$ can be obtained from $v$ by permuting finitely many entries of $v$. Consider the following binary relation:
\begin{align*}
 u  \succsim_{\text{SS}}  v \Longleftrightarrow  \text{  some finite permutation of $u$ dominates $v$. }  
 \end{align*}
This is the Suppes-Sen grading principle. It is the weakest preorder  satisfying {\bP} and {\bA} (see, e.g., \cite[p.~356]{BM07}).  We will use the fact that $\succsim_{\text{SS}}$ satisfies {\bS}, but violates {\bU}. 

For $c\in \mathbb{R}, n \in \mathbb{Z}_+=\{0, 1, 2, 3, \ldots\}$ and $u \in \U$, let 
\begin{align}
&([c]_n, u)=\begin{cases}
(\overbrace{c, c, \ldots, c}^{n \text{ times }}, u) & \text{ if $n\geq 1$},\\
  u & \text{ if $n=0$.}
\end{cases}
\end{align} 
Our last preorder is  defined as follows for all $u, v \in \U$:
\begin{align*}
  u \succsim_\ast v \Longleftrightarrow\text{    }  &([\bar{u}_\ast]_n, u) \succsim_{\text{SS}} ([\bar{v}^\ast]_m, v)   \text{  for some } n, m \in \mathbb{Z}_+.  
\end{align*}
Note that $\bar{u}_\ast$ and $\bar{u}^\ast$ (see  \eqref{def: limits}) are finite for each $u \in \U$.  

We first check that $\succsim_\ast$  is indeed a preorder. Reflexivity is obvious.   To show that $\succsim_\ast$ is transitive, let  $u, v, w \in \U$ be such that $u \succsim_\ast v$ and $v \succsim_\ast w$. That is, $([\bar{u}_\ast]_n, u) \succsim_{\text{SS}}  ([\bar{v}^\ast]_m, v)$ and $([\bar{v}_\ast]_k, v) \succsim_{\text{SS}}  ([\bar{w}^\ast]_l, w)$   for some $n, m, k, l \in \mathbb{Z}_+$. By the definition of $\succsim_{\text{SS}}$, we must then have that  $$\bar{u}^\ast \geq \bar{v}^\ast \geq \bar{w}^\ast \text{ and }\bar{u}_\ast \geq \bar{v}_\ast \geq \bar{w}_\ast.$$
 Since $\succsim_{\text{SS}}$ satisfies {\bS},  $([\bar{u}_\ast]_n, u) \succsim_{\text{SS}}  ([\bar{v}^\ast]_m, v)$ implies 
 \begin{align*}
 ( [\bar{v}_\ast]_{k}, [\bar{u}_\ast]_{n}, u) & \succsim_{\text{SS}}  \overbrace{([\bar{v}_\ast]_{k}, [\bar{v}^\ast]_{m},  v)}^{x}
 \intertext{ and  $([\bar{v}_\ast]_k, v) \succsim_{\text{SS}}  ([\bar{w}^\ast]_l, w)$  implies }
  \overbrace{([\bar{v}^\ast]_{m}, [\bar{v}_\ast]_{k}, v)}^{y} & \succsim_{\text{SS}}  ([\bar{v}^\ast]_{m},   [\bar{w}^\ast]_l, w).
 \end{align*}	
 By {\bA} and transitivity, $x \sim_{\text{SS}}  y$.  By transitivity,    \begin{align}\label{eq: klm} 
  ( [\bar{v}_\ast]_{k}, [\bar{u}_\ast]_{n}, u) \succsim_{\text{SS}}  ([\bar{v}^\ast]_{m},   [\bar{w}^\ast]_l, w).
  \end{align} 
    By {\bP}  and transitivity, \eqref{eq: klm} implies $([\bar{u}_\ast]_{n+k}, u) \succsim_{\text{SS}}   ([\bar{w}_\ast]_{l+m}, w)$, which means that $ u \succsim_\ast w$. Conclude that $\succsim_\ast$ is transitive and hence a preorder. 

To see that $\succsim_\ast$ satisfies {\bCP}, let $u \in \U$ be such that $\bar{u}$ is well defined, and let $v=(\bar{u}, u)$. Then $\bar{v}$ is well defined and  $\bar{u}=\bar{v}$. Since $([\bar{u}_\ast]_1, u) =([\bar{v}^\ast]_0, v)$, we have $([\bar{u}_\ast]_1, u) \succsim_{\text{SS}} ([\bar{v}^\ast]_0, v)$ and therefore $u \succsim_\ast v$.  Since $([\bar{v}_\ast]_0, v) =([\bar{u}^\ast]_1, u)$, we have $([\bar{v}_\ast]_0, v) \succsim_{\text{SS}} ([\bar{u}^\ast]_1, u)$ and therefore  $v \succsim_\ast u$.  Thus, $v=(\bar{u}, u) \sim_\ast u$, which shows that  $\succsim_\ast$ satisfies {\bCP}. 
 
To verify  {\bP}, suppose   $u$ is strictly better than $v$ by   {\bP}. Then  $([\bar{u}_\ast]_0, u) \succsim_{\text{SS}}  ([\bar{v}^\ast]_0, v)$, so $u \succsim_\ast v$.  To show that $v \neg \succsim_\ast u$,  we need to rule out the possibility that $([\bar{v}_\ast]_n,  v) \succsim_{\text{SS}} ([\bar{u}^\ast]_m, u)$ for some $n, m \in \mathbb{Z}_+$.  Suppose for contradiction that $([\bar{v}_\ast]_n,  v) \succsim_{\text{SS}} ([\bar{u}^\ast]_m, u)$ for  $n, m \in \mathbb{Z}_+$. Since $\succsim_{\text{AO}}$   satisfies  {\bP} and {\bA}, this implies that $([\bar{v}_\ast]_n,  v)  \succsim_{\text{AO}} ([\bar{u}^\ast]_m, u)$. By \eqref{eq CPs} and transitivity, we have $v  \succsim_{\text{AO}} ([\bar{v}_\ast]_n,  v)$ and  $([\bar{u}^\ast]_m, u) \succsim_{\text{AO}} u$. By transitivity, this means that $v  \succsim_{\text{AO}} u$, contradicting that $\succsim_{\text{AO}}$ satisfies {\bP}. We can therefore conclude that  $v \neg \succsim_\ast u$, so that  $u \succ_\ast v$. This shows that $\succsim_\ast$ satisfies {\bP}. 

It remains to show that $\succsim_\ast$ violates {\bU}. Define $u=(3, 5, 1, 1, 1, \ldots)$,  $v= (4, 2, 1, 1, 1, \ldots)$.  Since $u \succsim_{\text{SS}} v$, we have   $u \succsim_\ast v$. Let $\alpha=(-1, 1, 0, 0, 0, \ldots)$, let  $x= u+\alpha$, and let $y = v+\alpha$. Then $x=(2, 6, 1, 1,  1, \ldots)$ and $y=(3, 3, 1, 1, 1, \ldots)$.  Since $([\bar{x}_\ast]_n, x)=([1]_n, 2, 6, 1, 1,  1, \ldots)$ and  $([\bar{y}^\ast]_m, y)=([1]_m, 3, 3, 1, 1, 1, \ldots)$, there are no $n, m \in \mathbb{Z}_+$ with $([\bar{x}_\ast]_n, x) \succsim_{\text{SS}}  ([\bar{y}]^\ast_m, y)$.  Conclude that $x  \neg \succsim_\ast  y$. This shows that $\succsim_\ast$ violates  {\bU}.
 
The three preorders ($\succsim_{\neg \text{\bP}}$, $\succsim_\ast$ and $\succsim_{\text{O}}$) establish logical independence of   {\bP},  {\bU} and {\bCP}. These preorders also show that   Theorem \ref{theorem3} fails if any one of these three axioms is dropped.

\bibliographystyle{plain}
\bibliography{references-MDP}

\end{document}